\documentclass{amsart}
\usepackage{graphicx}
\usepackage{amscd}
\usepackage{amsmath}
\usepackage{amsfonts}
\usepackage{amssymb}
\theoremstyle{plain}
\newtheorem{theorem}{Theorem}
\newtheorem{corollary}[theorem]{Corollary}
\newtheorem{lemma}[theorem]{Lemma}
\newtheorem{proposition}[theorem]{Proposition}
\theoremstyle{definition}
\newtheorem{example}[theorem]{Example}

\newtheorem{definition}[theorem]{Definition}

\theoremstyle{remark}

\begin{document}
\title{A Schreier domain type condition}

\author{Zaheer Ahmad, Tiberiu Dumitrescu and Mihai Epure}
\address{Abdus Salam School of
Mathematical Sciences,
GC University Lahore, 68-B New
Muslim Town, Lahore 54600, Pakistan
}\email{zaheeir@yahoo.com}

\address{Facultatea de Matematica si Informatica,
University of Bucharest,
14 A\-ca\-de\-mi\-ei Str., Bucharest, RO 010014,
Romania}
\email{tiberiu@fmi.unibuc.ro, tiberiu\_dumitrescu2003@yahoo.com}

\address{Simion Stoilow Institute of Mathematics of the Romanian Academy
Research unit 5, P. O. Box 1-764, RO-014700 Bucharest, Romania}
\email{epuremihai@yahoo.com}


\begin{abstract}
\noindent 
We study the integral domains $D$  satisfying the following condition: whenever $I \supseteq AB$ with $I$, $A$, $B$  nonzero ideals, there exist  ideals $A'\supseteq A$ and $B'\supseteq B$ such that $I=A'B'$.
\end{abstract}

\maketitle
%
%

\noindent
In \cite{C}, Cohn introduced  the  notion of Schreier domain.  A domain $D$ is said to be a {\em Schreier domain} if $(1)$ $D$ is integrally closed and $(2)$  whenever $I, J_1, J_2$ are principal ideals of $D$ and $I\supseteq J_1J_2$, then  $I = I_1I_2$ for some principal ideals $I_1,I_2$ of $D$ with $I_i \supseteq J_i$ for $i = 1,2$.  The study of Schreier domains was continued in \cite{MR} and \cite{Zps} (where a domain was called  a {\em pre-Schreier domain} if it satisfies  condition $(2)$ above). In \cite{DM} and \cite{ADZ}, an  extension of the class of pre-Schreier domains was studied. A domain $D$ was called a {\em quasi-Schreier domain}  if whenever $I, J_1, J_2$ are invertible ideals of $D$ and $I\supseteq J_1J_2$, then  $I = I_1I_2$ for some (invertible) ideals $I_1,I_2$ of $D$ with $I_i \supseteq J_i$ for $i = 1,2$. 

In this paper we study the domains satisfying a Schreier-like condition for all nonzero  ideals. Since this class of domains turns out to be rather narrow, we use an ad hoc name for it.

\begin{definition}\label{41}
We call a domain $D$ a {\em  sharp domain} if whenever $I \supseteq AB$ with $I$, $A$, $B$  nonzero ideals of $D$, there exist  ideals $A'\supseteq A$ and $B'\supseteq B$ such that $I=A'B'$.
\end{definition}

If the domain $D$ is  Noetherian or Krull, then $D$ is sharp if and only if $D$ is a  Dedekind domain (Corollaries \ref{1} and \ref{42}).
In Proposition \ref{2}, we show that a sharp domain is  pseudo-Dedekind. In particular, a sharp domain is a completely integrally closed GGCD domain. The ring $E$ of entire functions is pseudo-Dedekind but not sharp (Example \ref{50}).
Recall (cf. \cite{Z} and \cite{AK}) that a domain $D$ is called a {\em pseudo-Dedekind domain} (the name used in \cite{Z}  was {\em generalized Dedekind domain}) if the $v$-closure of each nonzero ideal of $D$ is invertible. Also, recall from \cite{AA} that a  domain $D$ is called a {\em generalized GCD domain (GGCD domain)} if  the $v$-closure of each nonzero finitely generated ideal of $D$ is invertible. The definition of the $v$-closure is recalled below.
In Proposition \ref{5}, we show that a  valuation domain 
is sharp if and only if the value group  of $D$ is  a complete subgroup of the reals. 

The main results of this paper are Theorems \ref{9} and \ref{21}.
In Theorem \ref{9},  we show that the localizations of a  sharp domain at the maximal ideals are valuation domains with value group a complete subgroup of the reals. In par\-ti\-cu\-lar, a sharp domain is a Pr\"ufer   domain of dimension $\leq 1$. A key point in proving 
Theorem \ref{9} is the fact that if $D$ is a sharp domain and $x,y\in D-\{0\}$  such that $xD\cap yD=xyD$, then $xD+yD=D$ (Proposition \ref{7}).
The converse of Theorem \ref{9} is not true (Example \ref{44}).
In Theorem \ref{21}, we prove  the converse of Theorem \ref{9}  for the domains of finite character (i.e.,  domains whose every nonzero element is contained in only finitely many maximal ideals). 
The problem whether a sharp domain is necessarily of finite character is left open. A countable sharp domain is a Dedekind domain (Corollary \ref{45}).

For reader's convenience, we recall  the following facts. Let $D$ be a domain with quotient field $K$ and  $I$  a nonzero fractional ideal of $D$. The {\em $v$-closure}  of $I$ is the fractional ideal $I_v=(I^{-1})^{-1}$, where $I^{-1} = \{x \in K|\ xI \subseteq D\}$, and  $I$ is called a {\em $v$-ideal} if $I=I_v$.
The {\em $t$-closure}  of $I$ is the fractional ideal $I_t$ which is the union of the $v$-closures of the finitely generated nonzero 
subideals of $I$. Moreover,  $I$ is called a {\em $t$-ideal} if $I=I_t$. In general, we have $I\subseteq I_t\subseteq I_v$. A nonzero prime ideal $P$ of $D$ is called {\em $t$-prime} if $P=P_t$. For basic facts and  terminology not recalled in this paper, our references are \cite{G} and \cite{H}.
Throughout this paper, all rings are domains, that is, commutative, unitary and without zero-divisors.
\\[2mm]


We begin with a characterization of the sharp domains. If $I,H$ are ideals of a domain $D$, we denote by $I:H$ the ideal $\{x \in D \mid xH\subseteq I\}.$

\begin{proposition} \label{1}
A domain $D$ is sharp if and only if  for every two nonzero ideals $I,H$ we have $I=[I:(I:H)](I:H)$. 
\end{proposition}
{\em Proof.}   
$(\Rightarrow).$  Let $I,H$ be nonzero ideals  of $D$. Set  $A=I:(I:H)$ and $B=I:H$. Note that $AB\subseteq I$, $A=I:B$ and $I:A=I:(I:(I:H))=I:H=B$. 
As $D$ is sharp, there exists a 
factorization $I=A'B'$ with $A'$,$B'$ ideals such that $A'\supseteq A$ and $B'\supseteq B$. Then  $A\subseteq A'\subseteq I:B'\subseteq I:B=A$, so $A=A'$. Similarly, we get $B=B'$.
$(\Leftarrow).$ Let $I,A,B$ be nonzero ideals of $D$ such that $AB\subseteq I$. By our assumption, we get $I=[I:(I:A)](I:A)$. Note that $A\subseteq I:(I:A)$ and $B\subseteq I:A$. 
$\bullet$

\begin{corollary} \label{4}
A Dedekind domain is a sharp domain. 
\end{corollary}
{\em Proof.}
Let $D$ be a  Dedekind domain  and   $I,H$ nonzero ideals  of $D$.
Since $I:H$ is an invertible ideal, it follows easily that $I:(I:H)=I(I:H)^{-1}$. Hence $[I:(I:H)](I:H)=I(I:H)^{-1}(I:H)=I$. Apply Proposition \ref{1}.
$\bullet$\\[2mm]


\begin{proposition} \label{2}
Every sharp domain is  pseudo-Dedekind. In particular, a sharp domain is a completely integrally closed GGCD domain. 
\end{proposition}
{\em Proof.}
Let $D$ be a sharp domain,  $A$  a nonzero ideal of $D$ and $0\neq b \in A$. By  Proposition \ref{1}, $bD=[bD:(bD:A)](bD:A)$.  It follows that $bD:A=bA^{-1}$ is an invertible ideal, so $A^{-1}$ and $A_v$ are invertible ideals. Thus $D$ is pseudo-Dedekind. By \cite[Corollaries 1.4 and 1.5]{Z}, a pseudo-Dedekind domain is a completely integrally closed GGCD domain.
$\bullet$\\[2mm]

We show that for a pseudo-Dedekind domain $D$ it suffices to test the condition in Definition \ref{41} only for ideals $I$ with $I_v=D$. 

\begin{proposition} \label{3}
A pseudo-Dedekind domain $D$ is  sharp   if and only if
for all nonzero ideals  $I$,$A$,$B$   of $D$ such that $I\supseteq  AB$ and $I_v=D$, there exist  ideals  $A'\supseteq  A$ and $B'\supseteq  B$ such that $I=A'B'$.
\end{proposition}
{\em Proof.}   We prove the nontrivial implication. 
Let $I,A,B$ be nonzero ideals of $D$ such that $I\supseteq AB$. Then $I_v \supseteq A_vB_v$ and $I_v,A_v,B_v$ are invertible ideals, because $D$ is pseudo-Dedekind. 
A pseudo-Dedekind domain is a GGCD domain, cf. \cite[Corollary 1.5]{Z}, and a GGCD domain is quasi-Schreier, cf. \cite[Proposition 2.3]{DM}.
So there exist  invertible ideals  $A_1\supseteq  A_v$ and $B_1\supseteq  B_v$ such that $I_v=A_1B_1$. We have $I^{-1}=A_1^{-1}B_1^{-1}$, so $II^{-1} \supseteq (AA_1^{-1})(BB_1^{-1})$ and $AA_1^{-1}$, $BB_1^{-1}$ are integral ideals. Since $I_v$ is invertible, $(II^{-1})_v=D$. 
By our hypothesis, there exist  ideals  $A_2\supseteq  AA_1^{-1}$ and $B_2\supseteq  BB_1^{-1}$ such that $II^{-1}=A_2B_2$. Hence $I=(A_1A_2)(B_1B_2)$ and  $A_1A_2\supseteq  A$, $B_1B_2\supseteq  B$.
$\bullet$\\[2mm]

Next, we characterize the sharp valuation domains. 
Recall \cite[Exercise 21, page 551]{B}, that a {\em pseudo-principal domain} is a domain whose $v$-ideals are principal. Clearly, a quasi-local domain is pseudo-Dedekind if and only if it is pseudo-principal.

\begin{proposition} \label{5}
For a valuation domain $D$, the following assertions are equivalent:

$(a)$ $D$ is sharp.

$(b)$ $D$ is pseudo-Dedekind.

$(c)$ the value group  of $D$ is  a complete subgroup of the reals. 
\\
In particular, a sharp valuation domain has dimension $\leq 1$.
\end{proposition}
{\em Proof.}
$(b)\Leftrightarrow (c)$ is given in \cite{AK} at the bottom of pages 325 and 327 and 
$(a)\Rightarrow (b)$ follows from Proposition \ref{2}. We prove that 
$(b)$ and $(c)$ imply $(a)$. By Corollary \ref{4}, we may assume that the value group  of $D$ is the whole group of  real numbers. 
By Proposition \ref{3}, $D$ is sharp, because the maximal ideal is the only   proper ideal of $D$ whose $v$-closure is $D$. The ``in particular'' assertion follows from the well-known fact that a valuation domain has dimension $\leq 1$ if and only if its value group   is  a  subgroup of the reals (see \cite[page 45]{ZS}).
$\bullet$

\begin{proposition} \label{6}
If $D$ is a sharp domain, then every fraction ring $D_S$ of $D$ is also a sharp domain.
\end{proposition}
{\em Proof.} 
Let $I,A,B$ be nonzero ideals  of $D$ such that $ID_S\supseteq ABD_S$. Then $H=ID_S \cap D \supseteq  AB$. As $D$ is sharp, we get $H=A'B'$ with $A',B'$ ideals of $D$ such that  $A'\supseteq A$ and $B'\supseteq B$. Then $ID_S=HD_S=A'B'D_S$. 
$\bullet$

\begin{example}\label{50}
The ring $E$ of entire functions is pseudo-Dedekind but some localization of $E$ is not pseudo-Dedekind, cf.  \cite[Example 2.1]{Z}. By Proposition \ref{6}, $E$ is not a sharp domain.
\end{example}

\begin{proposition} \label{8}
If $D$ is a  sharp domain and $P$ is a $t$-prime ideal of $D$, then $D_P$ is a valuation domain whose value group   is  a complete subgroup of the reals. In particular, in a sharp domain every $t$-prime ideal of $D$ has height one.
\end{proposition}
{\em Proof.} 
By Proposition \ref{2}, $D$ is a GGCD domain.  By \cite[page 218]{AA}, \cite[Corollary 4.3]{MZ} and Proposition \ref{6}, $D_P$ is a sharp valuation domain. Apply Proposition \ref{5}.
$\bullet$\\[2mm]

Recall that two nonzero  elements $x,y$ of a domain $D$ are called {\em $v$-coprime} if $(xD+yD)_v=D$ (equivalently $xD\cap yD=xyD$, equivalently   $xD : yD=xD$).

\begin{proposition} \label{7}
Let $D$ be a sharp domain and $x,y$ two nonzero $v$-coprime elements. Then $xD+yD=D$.
\end{proposition}
{\em Proof.}  
We have $(x,y)^2 \subseteq (x^2,y)$, so $(x^2,y)=AB$ with $A,B$ ideals such that $A,B\supseteq (x,y)$. Note that $(x^2,y):(x,y)=(x,y)$. Indeed, if $a \in (x^2,y):(x,y)$, then $ax=bx^2+cy$ for some $b,c \in D$, so $c \in xD : yD=xD$, hence $a=bx+(c/x)y$ belongs to $(x,y)$.
From $(x^2,y)=AB$, we get $A\subseteq (x^2,y):B\subseteq (x^2,y):(x,y)=(x,y)$, so $A=(x,y)$. Similarly, we get $B=(x,y)$. Then $(x^2,y)=(x,y)^2$. So $y=fx+gy^2$ for some $f,g \in D$, hence  $f \in yD : xD=yD$, thus $1=(f/y)x+gy$, that is, $xD+yD=D$.
$\bullet$\\[2mm]

Theorems \ref{9} and \ref{21} are the main results of this paper.

\begin{theorem}\label{9}
If $D$ is a sharp domain, then   $D_M$ is a  valuation domain with value group a complete subgroup of the reals, for each maximal ideal $M$ of $D$.
In par\-ti\-cu\-lar, a sharp domain is a Pr\"ufer   domain of dimension $\leq 1$.
\end{theorem}
{\em Proof.}
By Proposition \ref{6}, we may assume that $D$ is  quasi-local  with nonzero maximal ideal $M$.
Suppose that the height of $M$ is $\geq 2$.  By Proposition 
\ref{2}, $D$ is a quasi-local GGCD  domain, hence a GCD domain, cf. \cite[Corollary 1]{AA}.
By Proposition \ref{8}, $M$ is not a $t$-ideal, so 
$M_t=D$. Since $D$ is a GCD domain, there exist two $v$-coprime elements $x,y \in M$ (see the paragraph before Theorem 4.8 in \cite{A}). But this  contradicts  Proposition \ref{7}. It remains that $M$ has height one, hence it is a $t$-prime, cf. \cite[Proposition 6.6]{H}. Now apply Proposition \ref{8} to conclude. The ``in particular'' assertion is clear.
$\bullet$\\[2mm]

According to \cite{HZ}, a {\em TV domain} is  a domain in which every $t$-ideal is a $v$-ideal. It is well known that  Noetherian domains and  Krull  domains are TV domains.

\begin{corollary}\label{42}
If $D$ is a sharp TV domain, then   $D$ is a  Dedekind domain. In particular, if a  sharp domain is Noetherian  or Krull, then it is a Dedekind domain.
\end{corollary}
{\em Proof.}
Let $D$ be a sharp TV domain. By Theorem \ref{9}, $D$ is a Pr\"ufer   domain, so every nonzero ideal of $D$ is a $t$-ideal, hence a $v$-ideal, because $D$ is a TV domain. Since  $D$  is also a   pseudo-Dedekind domain (cf. Proposition  \ref{2}), it follows that every nonzero ideal of $D$ is invertible. Thus $D$ is a  Dedekind domain.
$\bullet$\\[2mm]

The converse of Theorem \ref{9} is not true.
Recall \cite[page 434]{G} that  a domain $D$  is said to be  {\em almost Dedekind} if $D_M$ is a discrete (Noetherian) valuation domain for each maximal ideal $M$ of $D$. 
We exhibit an almost Dedekind domain which is not  a sharp domain (not even pseudo-Dedekind).

\begin{example}\label{44}
Let $D$ be the almost Dedekind domain constructed in the proof of \cite[Proposition 7]{CD}. We recall some properties of $D$ proved there. The maximal ideals of $D$ are the principal ideals $(p_iD)_{i\geq 1}$ and the ideal $M=(q_0,q_1,...,q_n,...)$. Here  $(q_i)_{i\geq 0}$ are nonzero elements of $D$ such that $q_{i-1}=p_iq_i$ and $p_i$ does not divide $q_i$ for all $i\geq 1$. 
Note that $M$  is not finitely generated, because it is the union of the strictly ascending chain of principal ideals  $(q_iD)_{i\geq 0}$.
We claim that $D$ is not pseudo-Dedekind, so it is not a sharp domain (cf. Proposition \ref{3}). For that, it suffices to prove that the $v$-ideal $\cap_{i\geq 1}p_{2i-1}D$  equals the union of the strictly ascending chain of principal ideals  $p_1q_2D\subset$ $p_1p_3q_4D\subset p_1p_3p_5q_6D\cdots$, so it is not finitely generated.
Indeed, the inclusion $\supseteq$ is clear. Conversely, let $x\in \cap_{i\geq 1}p_{2i-1}D$. If $x\notin M$, then $1=ax+bq_{2n}$ for some $a,b\in D$ and $n\geq 0$. But this is a contradiction, because  $p_{2n+1}$ divides both $x$ and $q_{2n}$. So $x\in M$, say $x=cq_{2n}$ for some $c\in D$ and $n\geq 1$. Since $x\in \cap_{i\geq 1}p_{2i-1}D$ and $q_{2n}$ is not divisible by $p_1$,$p_3$,...,$p_{2n-1}$, we get that $x\in p_1p_3\cdots p_{2n-1}q_{2n}D$.
\end{example}

We give a partial converse of Theorem \ref{9}. Recall  that 
 a domain $D$ is said to be of {\em finite character} if every nonzero element is contained in only finitely many maximal ideals.
And $D$ is said to be  {\em h-local} if $D$ is of finite character and every nonzero prime ideal of $D$  is contained in  a unique maximal ideal of $D$. It is easy to see  that a one-dimensional domain of finite character is h-local. The next lemma was implicit  in \cite[Proposition 3.1]{O}.

\begin{lemma} \label{10}
Let  $D$ be a h-local domain, $A,B$ nonzero ideals of $D$ and $M\in Max(D)$. Then $(A:B)D_M=AD_M:BD_M$.
\end{lemma}
{\em Proof.} 
Let $K$ denote the quotient field of $D$.
The inclusion $(\subseteq)$ is clear. Conversely, let $x\in AD_M:BD_M$. We may assume that $x\in D$. Pick $a\in A-\{0\}$. Since $D$ is h-local, we have $[M]D_M=K$ where $[M]=\cap \{ D_N\mid N\in Max(D)$ and $N\neq M\}$, cf. \cite[Proposition 3.1]{O}. Consequently, 
there exist $y\in [M]$ and $s\in D-M$ such that $x/a=y/s$. So $sx=ay$. Note that $ayB\subseteq AD_N$ for each $N\in Max(D)-\{M\}$. So $sxB_Q\subseteq AD_Q$ for each $Q\in Max(D)$, hence $sx\in A:B$. Thus $x\in (A:B)D_M$.
$\bullet$\\[2mm]

We show that the converse of Theorem \ref{9} is true for  a domain  of finite character.

\begin{theorem} \label{21}
Let  $D$ be a  domain of finite character such that $D_M$ is a  valuation domain with value group a complete subgroup of the reals  for each $M\in Max(D)$. Then $D$ is a sharp domain.
\end{theorem}
{\em Proof.} 
Let $I,A$ be nonzero ideals of $D$. 
By Proposition \ref{1}, it suffices to check locally that $(I:A)[I:(I:A)]=I$. Let $M$ be a maximal ideal of $D$. 
Since $D$ is one-dimensional  of finite character, it is h-local.
By Lemma \ref{10}, we have $(I:A)[I:(I:A)]D_M=(ID_M:AD_M)[ID_M:(ID_M:AD_M)]=ID_M$, where the last equality follows from Propositions \ref{5} and \ref{1}. 
$\bullet$\\[2mm]

We do not know if a sharp domain is necessarily of finite character. A connected question, which is up to our knowleadge not solved, is whether a pseudo-Dedekind almost Dedekind domain is necessarily a Dedekind domain. We end our paper with two results for countable domains.

\begin{proposition}\label{22}
If $D$ is a countable pseudo-Dedekind Pr\"ufer domain, then $D$ is of finite character.
\end{proposition}
{\em Proof.} 
Assume that $D$ is not of finite character. By \cite[Corollary 7]{DZ},  there exists a nonzero element $z$ and an infinite family $(I_n)_{n\geq 1}$ of invertible proper mutually comaximal ideals containing $z$. For each nonempty set of natural numbers $\Lambda$, consider the $v$-ideal $I_\Lambda=\cap_{n\in \Lambda} I_n$ (note that 
$I_\Lambda$ contains $z$). As $D$ is  pseudo-Dedekind, $I_\Lambda$ is invertible. We claim that $I_\Lambda\neq I_{\Lambda'}$ whenever $\Lambda$, $\Lambda'$ are distinct nonempty sets of natural numbers.
Deny. Then there exists a nonempty set  of natural numbers $\Gamma$ and some $k\notin \Gamma$ such that $I_k\supseteq I_\Gamma$. Consider the ideal  $H=I_k^{-1} I_\Gamma\supseteq I_\Gamma$. If $n\in \Gamma$, then $I_n \supseteq I_\Gamma =I_kH$, so $I_n \supseteq H$, because $I_n+I_k=D$. It follows that $I_\Gamma\supseteq H$, so $I_\Gamma = H=I_k^{-1} I_\Gamma$. Since $I_\Gamma$ is invertible, we get $I_k=D$, a contradiction. Thus the claim is proved. But then it follows that $\{I_\Lambda\mid \emptyset\neq \Lambda\subseteq \mathbb{N}\}$ is an uncountable set of invertible ideals. This leads to a contradiction, because $D$ being countable, it has countably many finitely generated ideals.
$\bullet$

\begin{corollary}\label{45}
If $D$ is a countable sharp domain, then $D$ is a Dedekind domain.
\end{corollary}
{\em Proof.} 
We may assume that $D$ is not a field.
By Theorem \ref{9}, $D$ is a Pr\"ufer   domain. Now Propositions \ref{2} and \ref{22} show that $D$ is of finite character. Let $M$ be a maximal ideal of $D$. By Theorem \ref{9}, $D_M$ is a countable valuation domain with value group $\mathbb{Z}$ or $\mathbb{R}$, so $D_M$ is a DVR. Thus $D$ is a Dedekind domain.
$\bullet$\\[2mm]

{\bf Acknowledgements.} The first author was partially supported by an HEC (Higher Education Commission, Pakistan)  grant. The second author gratefully acknowledges the warm 
hospitality of the Abdus Salam School of Mathematical Sciences GCU Lahore during his visits in 2006-2010. The third author was supported by UEFISCDI, project number 83/2010, PNII-RU code TE\_46/2010, program Human Resources.
\\[7mm]

\end{document}